# A Comment on the Conformable Euler's Finite Difference Method

D.P. Clemence-Mkhope

Department of Mathematics and Statistics, North Carolina A&T State University, Greensboro, NC 27411

E-mail: clemence@ncat.edu

**Abstract**

A method, recently advanced as the conformable Euler's method, a general method for the finite difference discretization of fractional initial value problems for $0 < \alpha \leq 1$, is shown to be valid only for $\alpha = 1$. The property of the conformable fractional derivative used to show this limitation of the method is used, together with the integer definition of the derivative, to derive a modified conformable Euler's method for the initial value problem considered.

**Keywords:** conformable fractional derivative, conformable Euler's method, fractional finite difference, fractional initial value problem

## 1. Introduction

The conformable fractional derivative (CFD) was introduced by Khalil et al. [1] and, in a short time, its theory has been developed by various authors and the definition widely applied (see references in [2] for example). The CFD is defined as

$$_0^C T_t^\alpha (f)(t) = \lim_{\varepsilon \to 0} \frac{f(t+\varepsilon t^{1-\alpha}) - f(t)}{\varepsilon}.$$

It exists if and only if $\frac{d}{dt} f(t)$ exists, and satisfies the identity

$$_0^C T_t^\alpha (f)(t) = t^{1-\alpha} \frac{d}{dt} f(t). \tag{1.1}$$

Since most conformable fractional differential equations systems do not have exact analytic solutions, numerical approximation methods must be developed. Among the many numerical methods that have been proposed to solve the conformable fractional differential equations, the finite difference method

$$y_{k+1} - y_k = \frac{1}{\alpha} h^\alpha f(t_k, y(t_k)), \ 0 \leq t \leq N, \text{ where } h = \frac{b-a}{N}, \tag{1.2}$$

has been used recently to discretize the fractional initial value problem:

$$D_t^\alpha y(t) = f(t; y(t)), y(t_0) = y_0, \ a \leq t \leq b, \tag{1.3}$$

and some generalizations. Termed the conformable Euler's method, it is casually adopted in [2] and justified by applying the fractional power series expansion. In [2], the method is used to solve conformable fractional differential equations system with time delays, and in [3], it is used in conjunction with a to calculate numerical solutions for testing the hyperchaos of conformable derivative models for financial systems.

The purpose of this short note is to show that the method (1.2) as a model for (1.3) is valid only for $\alpha = 1$, and to propose a modified conformable Euler's method. Both are accomplished using the expression of the CFD in terms of the ordinary, $\alpha = 1$, derivative.

## 2. The Conformable Euler's Method

### 2.1 Justification of the Conformable Euler's Method

The method (1.2), referred to in [2] and [3] as the conformable Euler's method for (1.3), is obtained by considering a power series expansion as follows. Since $h = t_{k+1} - t_k$, it is assumed that there exist $\theta_k$ where $0 < \theta_k < 1$ is such that

$$y(t_{k+1}) - y(t_k) = \frac{1}{\alpha} h^\alpha (D_t^\alpha y)(t_k) + \frac{1}{2\alpha^2} h^{2\alpha} (D_t^{2\alpha} y)(t_k + \theta_k h). \tag{2.1}$$

Letting $y(t_{k+1}) - y(t_k) \to y_{k+1} - y_k$ and substituting $(D_t^\alpha y)(t_k) = f(t_k, y_k)$ into (2.1) results in

$$y_{k+1} - y_k = \frac{1}{\alpha} h^\alpha f(t_k, y_k) + \frac{1}{2\alpha^2} h^{2\alpha} (D_t^{2\alpha} y)(t_k + \theta_k h),$$

or, equivalently

$$\alpha \frac{y_{k+1} - y_k}{h^\alpha} = f(t_k, y_k) + \frac{1}{2\alpha} h^\alpha (D_t^{2\alpha} y)(t_k + \theta h). \tag{2.2}$$

For $h$ small enough, ignoring the second term on the right-hand side of (2.2) yields the conformable Euler's method (1.2):

$$y_{k+1} - y_k = \frac{1}{\alpha} h^\alpha f(t_k, y_k),$$

which reduces to the usual Euler's method for $\alpha = 1$.

### 2.2 Validity of the Conformable Euler's Method

Since the CFD satisfies property (1.1), substituting $(D_t^\alpha y)(t_k) = (t_k)^{1-\alpha} \frac{dy}{dt}(t_k)$ into (2.1) results in

$$y_{k+1} - y_k = \frac{1}{\alpha} h^\alpha (t_k)^{1-\alpha} \frac{dy}{dt}(t_k) + \frac{1}{2\alpha^2} h^{2\alpha} (D_t^{2\alpha} y)(t_k + \theta_k h),$$

or, equivalently

$$\alpha \frac{y_{k+1}-y_k}{h^\alpha} = (t_k)^{1-\alpha} \frac{dy}{dt}(t_k) + \frac{1}{2\alpha} h^\alpha (D_t^{2\alpha} y)(t_k + \theta h). \tag{2.3}$$

For $h$ small enough, ignoring the second term on the right-hand side of (2.3) yields

$$\alpha \frac{y_{k+1}-y_k}{h^\alpha} = (t_k)^{1-\alpha} \frac{dy}{dt}(t_k),$$

and therefore

$$\frac{\alpha}{h^{\alpha-1}} \frac{y_{k+1}-y_k}{h} = (t_k)^{1-\alpha} \frac{y_{k+1}-y_k}{h},$$

from which it follows that

$$\alpha h^{1-\alpha} = (t_k)^{1-\alpha}. \tag{2.4}$$

Next, let us consider separately the cases of (2.4) (a) $t_0 = 0$ and (b) $t_0 \neq 0$

(a) If $t_0 = 0$, then $t_k = kh$, so that

$$\alpha h^{1-\alpha} = (t_k)^{1-\alpha} = (kh)^{1-\alpha} = k^{1-\alpha} h^{1-\alpha},$$

from which we conclude that

$$\alpha = k^{1-\alpha},$$

which is true for all $k$ only for $\alpha = 1$.

(b) If $t_0 \neq 0$, then $t_k = t_0 + kh$, so that

$$\alpha h^{1-\alpha} = (t_k)^{1-\alpha} = (t_0 + kh)^{1-\alpha} = h^{1-\alpha} \left(\frac{t_0}{h} + k\right)^{1-\alpha},$$

from which we conclude that

$$\alpha = \left(\frac{t_0}{h} + k\right)^{1-\alpha}, \tag{2.5}$$

which is true for all $h > 0, t_0, k$ only for $\alpha = 1$. Note that, if $\alpha < 1$ is assumed on the right-hand side of (2.5), then for fixed $t_0, k$, and writing the left-hand side as $\alpha(h)$, there results that $\alpha(h) \to \infty$ as $h \to 0$, a contradiction.

Since (2.4) holds if and only if both (1.1) and (1.2) hold, we conclude therefore from (a) and (b) that both (1.1) and (1.2) hold if and only if $\alpha = 1$.

### 3. The Modified Conformable Euler's Method

To obtain a method consistent with property (1.1), rewrite $t^{1-\alpha}$ as a derivative and then use the $\alpha = 1$ definition of the derivative to get

$${}_0^C T_t^\alpha y(t) = t^{1-\alpha} \frac{d}{dt} y(t) = \left(\frac{d}{dt} y(t)\right) / \left(\frac{1}{\alpha} \frac{d}{dt}(t^\alpha)\right) = \lim_{h \to 0} \left(\frac{y(t+h)-y(t)}{h}\right) / \left(\frac{1}{\alpha} \frac{(t+h)^\alpha - t^\alpha}{h}\right)$$

$$= \alpha \lim_{h \to 0} \left(\frac{y(t+h)-y(t)}{(t+h)^\alpha - t^\alpha}\right) \tag{3.1}$$

For $h$ small enough, therefore, and making the identifications

$$t \to t_k, t + h \to t_{k+1}, y(t+h) \to y_{k+1}, y(t) \to y_k$$

in (3.1) result in the following discrete representation of ${}_0^C T_t^\alpha y(t)$:

$${}_0^C T_t^\alpha y(t) = \alpha \lim_{h \to 0} \left( \frac{y(t+h) - y(t)}{(t+h)^\alpha - t^\alpha} \right) \to \alpha \frac{y_{k+1} - y_k}{(t_{k+1})^\alpha - (t_k)^\alpha}.$$

and the modified conformable Euler's method for (1.3) is therefore given by

$$y_{k+1} - y_k = \frac{1}{\alpha}[(t_{k+1})^\alpha - (t_k)^\alpha] f(t_k, y_k),$$

valid for $0 < \alpha \leq 1$, which is also a generalization of the Euler method for $\alpha = 1$.

**Conclusion**

A discretization method for the fractional initial value problem, with the conformable fractional derivative, has been considered that extends the integer Euler method and is termed the conformable Euler's method in recent literature; its justification using a fractional series expansion is recalled. It has been shown that the assumption ${}_0^C T_t^\alpha (f)(t) = t^{1-\alpha} \frac{d}{dt} f(t)$ leads to the conclusion that the method is valid only for $\alpha = 1$. A modified conformable Euler's method is proposed that is derived from rewriting the term of the right-hand side of the assumption above as a quotient of ordinary derivatives and then using the integer definition of the derivative.

**Funding:**     This research did not receive any specific grant from funding agencies in the public, commercial, or not-for-profit sectors.